\documentclass[12pt]{article}

\usepackage{amscd,amsmath, amssymb, fancyhdr,
  epsfig,color}
\definecolor{dblue}{rgb}{0,0,.6}

\usepackage[small,nohug,heads=littlevee]{diagrams}
\diagramstyle[labelstyle=\scriptstyle]

\usepackage[backref=page]{hyperref}
\renewcommand*{\backref}[1]{}
\renewcommand*{\backrefalt}[4]{%
    \ifcase #1 (Not cited.)%
    \or        (Cited on page~#2.)%
    \else      (Cited on pages~#2.)%
    \fi}

\hypersetup{
     colorlinks   = true,
     citecolor    = magenta,
     linkcolor= blue
}

\numberwithin{equation}{section}

\newcommand{\version}{version 1.0,\ \   Aug. 14, 2019}

\def\eqref#1{(\ref{#1})}
\newcommand{\goth}{\mathfrak}

\newcommand{\arrow}{{\:\longrightarrow\:}}
\newcommand{\Z}{{\Bbb Z}}
\newcommand{\C}{{\Bbb C}}

\newcommand{\R}{{\Bbb R}}
\newcommand{\Q}{{\Bbb Q}}

\newcommand{\6}{\partial}
\def\1{\sqrt{-1}\:}
\newcommand{\restrict}[1]{{\left|_{{\phantom{|}\!\!}_{#1}}\right.}}
\newcommand{\cntrct}                
{\hspace{2pt}\raisebox{1pt}{\text{$\lrcorner$}}\hspace{2pt}}

\makeatletter
\def\x@arrow{\DOTSB\Relbar}
\def\xlongequalsignfill@{\arrowfill@\x@arrow\Relbar\x@arrow}
\newcommand{\xlongequal}[2][]{%
        \ext@arrow 0099\xlongequalsignfill@{#1}{#2}}
\def\xlongrightarrowfill@{\arrowfill@\relbar\relbar\longrightarrow}
\newcommand{\xlongrightarrow}[2][]{%
        \ext@arrow 0099\xlongrightarrowfill@{#1}{#2}}
\makeatother

\newcommand{\calo}{{\cal O}}

\renewcommand{\bar}{\overline}
\renewcommand{\phi}{\varphi}
\renewcommand{\epsilon}{\varepsilon}
\renewcommand{\geq}{\geqslant}
\renewcommand{\leq}{\leqslant}

\newcommand{\tors}{\operatorname{\text{\sf\small tors}}}

\newcommand{\Tot}{\operatorname{Tot}}
\newcommand{\Id}{\operatorname{Id}}

\newcommand{\const}{\operatorname{\text{\sf const}}}

\newcommand{\Per}{\operatorname{Per}}

\newcommand{\Hilb}{\operatorname{Hilb}}

\newcommand{\Gr}{\operatorname{Gr}}
\newcommand{\Def}{\operatorname{Def}}

\newcommand{\Mor}{\operatorname{Mor}}
\newcommand{\St}{\operatorname{St}}

\renewcommand{\Re}{{\operatorname{Re}}}
\renewcommand{\Im}{{\operatorname{Im}}}

\newcommand{\Ham}{\operatorname{Ham}}
\newcommand{\Perspace}{\operatorname{{\Bbb P}\sf er}}

\newcounter{Mycounter}[section]
\newcounter{lemma}[section]
\renewcommand{\thelemma}{{Lemma \thesection.\arabic{lemma}}}
\newcommand{\lemma}{%
    \setcounter{lemma}{\value{Mycounter}}
    \refstepcounter{lemma}
    \stepcounter{Mycounter}
    {\noindent \bf \thelemma:\ }}

\newcounter{claim}[section]
\renewcommand{\theclaim}{{Claim \thesection.\arabic{claim}}}
\newcommand{\claim}{%
    \setcounter{claim}{\value{Mycounter}}
    \refstepcounter{claim}
    \stepcounter{Mycounter}
    {\noindent \bf \theclaim:\ }}

\newcounter{sublemma}[section]
\newcounter{corollary}[section]
\renewcommand{\thecorollary}{{Corollary \thesection.\arabic{corollary}}}
\newcommand{\corollary}{%
    \setcounter{corollary}{\value{Mycounter}}
    \refstepcounter{corollary}
    \stepcounter{Mycounter}
    {\noindent \bf \thecorollary:\ }}

\newcounter{theorem}[section]
\renewcommand{\thetheorem}{{Theorem \thesection.\arabic{theorem}}}
\newcommand{\theorem}{%
    \setcounter{theorem}{\value{Mycounter}}
    \refstepcounter{theorem}
    \stepcounter{Mycounter}
    {\noindent \bf \thetheorem:\ }}

\newcounter{conjecture}[section]
\renewcommand{\theconjecture}{{Conjecture \thesection.\arabic{conjecture}}}
\newcommand{\conjecture}{%
    \setcounter{conjecture}{\value{Mycounter}}
    \refstepcounter{conjecture}
    \stepcounter{Mycounter}
    {\noindent \bf \theconjecture:\ }}

\newcounter{proposition}[section]
\renewcommand{\theproposition}
      {{Proposition \thesection.\arabic{proposition}}}
\newcommand{\proposition}{%
    \setcounter{proposition}{\value{Mycounter}}
    \refstepcounter{proposition}
    \stepcounter{Mycounter}
    {\noindent \bf \theproposition:\ }}

\newcounter{definition}[section]
\renewcommand{\thedefinition}
      {{Definition~\thesection.\arabic{definition}}}
\newcommand{\definition}{%
    \setcounter{definition}{\value{Mycounter}}
    \refstepcounter{definition}
    \stepcounter{Mycounter}
    {\noindent \bf \thedefinition:\ }}

\newcounter{example}[section]
\renewcommand{\theexample}{{Example \thesection.\arabic{example}}}
\newcommand{\example}{%
    \setcounter{example}{\value{Mycounter}}
    \refstepcounter{example}
    \stepcounter{Mycounter}
    {\noindent \bf \theexample:\ }}

\newcounter{remark}[section]
\renewcommand{\theremark}{{Remark \thesection.\arabic{remark}}}
\newcommand{\remark}{%
    \setcounter{remark}{\value{Mycounter}}
    \refstepcounter{remark}
    \stepcounter{Mycounter}
    {\noindent \bf \theremark:\ }}

\newcounter{problem}[section]
\newcounter{question}[section]
\newcommand{\proof}{{\bf Proof:\:}}

\makeatletter

\@addtoreset{equation}{section} \@addtoreset{footnote}{section} \makeatother

\def\blacksquare{\hbox{\vrule width 5pt height 5pt depth 0pt}}
\def\endproof{\blacksquare}

\begin{document}
\begin{center}
{\LARGE\bf
Deformations and BBF form\\[4mm] on non-K\"ahler  holomorphically\\[4mm] symplectic manifolds}
\medskip

Nikon Kurnosov\footnote{ Nikon Kurnosov is partially supported by
by  the  Russian Academic Excellence Project '5-100', foundation ``BASIS'', Simons Travel grant, and by the contest ``Young Russian Mathematics''.},
Misha Verbitsky\footnote{ Misha Verbitsky is partially supported by
by  the  Russian Academic Excellence Project '5-100', FAPERJ E-26/202.912/2018 
and CNPq - Process 313608/2017-2}

\end{center}

{\small \hspace{0.02\linewidth}
\begin{minipage}[t]{0.85\linewidth} \small
{\bf Abstract}
In 1995, Dan Guan constructed examples of non-K\"ahler, simply-connected
holomorphically symplectic manifolds. An alternative construction,
using the Hilbert scheme of Kodaira-Thurston surface, was given
by F. Bogomolov. We investigate topology and deformation theory
of Bogomolov-Guan manifolds and show that it is similar
to that of hyperk\"ahler manifolds. We prove the local Torelli theorem, showing
that holomorphically symplectic deformations of BG-manifolds
are unobstructed, and the corresponding period map is
locally a diffeomorphism. Using the local Torelli theorem,
we prove the Fujiki formula for a BG-manifold $M$, showing that
there exists a symmetric form $q$ on $H^2(M)$ such that
for any $w\in H^2(M)$ one has $\int_M w^{2n}=q(w,w)^n$.
This form is a non-K\"ahler version of  
the Beauville-Bogomolov-Fujiki form
known in hyperk\"ahler geometry.
\end{minipage}
}

{\small 
\tableofcontents
}


\section{Introduction}


\subsection{Deformation theory for holomorphically symplectic manifolds}

Bogomolov-Guan manifolds are simply-connected, compact,
holomorphically symplectic manifolds, obtained as follows.
Kodaira-Enriques classification implies that
any holomorphically symplectic complex surface is isomorphic
to a torus, K3 surface, or Kodaira-Thurston surface which is 
defined as follows. Let $L$ be a line bundle
on an elliptic curve $E$ with the first Chern class
$c_1(L)\neq 0$. Denote by $\tilde S$ the corresponding $\C^*$-bundle
on $E$ obtained by removing the zero section,
$\tilde S=\Tot(L) \backslash 0$. Fix a complex
number $\lambda$ with $|\lambda| > 1$, and
let $h_\lambda:\; \tilde S \arrow \tilde S$
be the corresponding homothety of $\tilde S$.
The quotient $\tilde S/\langle h_\lambda\rangle$
is called {\bf a primary Kodaira-Thurston surface};
in the present paper we shall abbreviate this
to {\bf Kodaira-Thurston surface}.
It is an isotrivial elliptic fibration over 
the elliptic curve $E$, with the fiber
identified with an elliptic curve $E_L:\;= \C^*/\langle \lambda\rangle$.
Therefore, it is holomorphically symplectic.

Using K3 surface and a compact torus surface, one can
construct two famous families of hyperk\"ahler manifolds,
due to A. Beauville (\cite{_Beauville_}).
The Hilbert scheme of points on K3 surface
(itself obtained by a resolution of singularities
of the symmetric power of K3) is the first. The second,
known as {\bf generalized Kummer variety}, is
obtained from the Hilbert scheme of a torus
by taking a two-sheeted cover of its Albanese fiber.

Bogomolov-Guan manifolds were obtained in a similar
way starting from a Kodaira-Thurston surface $S$.
Let $S^{[n]}$ be the Hilbert scheme of $S$, 
$S^{(n)}$ its symmetric power, and $\pi_S:\; S \arrow E$
the elliptic fibration constructed above. 
Applying $\pi_S$ to each component of $S^{(n)}$
and summing up, we obtain a holomorphic
projection from $S^{(n)}$ to $E$;
taking the composition with 
the resolution $r:\; S^{[n]} \arrow S^{(n)}$,
we obtain an isotrivial fibration
$\pi:\; S^{[n]}\arrow E$. Denote its fiber by $F^{[n]}$.
The Hilbert scheme of a holomorphically symplectic
surface is again holomorphically symplectic (\cite{_Beauville_}).
Therefore, $F^{[n]}$ is a smooth divisor in 
a holomorphically symplectic manifold $(S^{[n]}, \Omega)$.
The restriction of $\Omega$ to $F^{[n]}$ has rank 
$2n-2$, because $F^{[n]}\subset S^{[n]}$ is a divisor.

Let $E_L$ be the fiber of  the projection of $S\arrow E$.
Then elliptic curve $E_L$, considered as a group,
acts on $S$; we extend this action to the Hilbert scheme
$S^{[n]}$ in a natural way.

Denote by $K\subset TF^{[n]}$ the kernel of $\Omega\restrict {F^{[n]}}$,
that is, the set of all $x\in TF^{[n]}$ such that 
$\Omega\restrict {F^{[n]}}(x, \cdot)=0$. 
The corresponding foliation is called {\bf the characteristic foliation};
this foliation is defined on a divisor in any
holomorphically symplectic manifold. It is not hard to 
see that these leaves are identified with orbits of the
natural $E_L$-action on $S^{[n]}$.

The leaf space $W$ of $K$ 
is a smooth holomorphically symplectic 
orbifold, but it is never smooth.
However, under additional numerical
assumptions, $W$ has a smooth finite covering, ramified
in the singular points of $W$. Suppose that
the degree of the line bundle $L$ over $E$ is divisible by $n$.
Bogomolov has shown that in this case the leaf space $W$
admits an order $n^2$ cyclic covering which is
smooth, simply connected and holomorphically
symplectic. This covering is called 
{\bf the Bogomolov-Guan manifold};
for a more detailed exposition see 
Subsection \ref{BG-manifolds-subsection}.

Discovery of Bogomolov-Guan manifolds was unexpected and quite
surprising, because it was assumed that such manifold
cannot exist. Indeed, in an MPIM 
preprint \cite{_Todorov:MPIM_}, A. Todorov claimed
that a simply-connected holomorphically symplectic manifold
with $H^{2,0}(M)=\C$ is always hyperk\"ahler.

The deformation theory of compact Calabi-Yau manifolds of 
K\"ahler type is well understood since the pioneering works
of Bogomolov, Tian and Todorov (\cite{_Bogomolov:TT_,_Tian:TT_,_Todorov:TT_}).
Deformations of Calabi-Yau manifolds satisfy the local Torelli
theorem. Namely, these deformations are unobstructed, 
their Kuranishi space $\Def(M)$ is smooth, and the Kodaira-Spencer
map \[ KS:\; T_I \Def(M) \arrow H^1(TM)
\] is an isomorphism.

The proof of this result is based on Tian-Todorov lemma 
(\ref{_TT_Lemma_}), which is given by a local and fairly elementary
computation, and the $\6\bar\6$-lemma, which uses deep
arguments of K\"ahler geometry. 
As shown in \cite{_Iacono_}, the $\6\bar\6$-lemma
can be replaced by a more abstract and general condition,
the homotopy abelianness of the differential graded
algebra $(\Lambda^{0,*}(TM), \bar\6)$ which is responsible
for the deformations. However, without the K\"ahler
condition, the local Torelli theorem is known to fail
(\cite{_Ghys_}). 

An alternative approach to deformations
of holomorphically symplectic manifolds was
suggested in \cite{_KV:deformations_}.
Instead of studying the deformations of
complex structures, \cite{_KV:deformations_} uses
holomorphically symplectic deformations, that is,
deformations of the holomorphically symplectic
form (the complex structure is determined
by the holomorphically symplectic form uniquely).
The Kuranishi deformation space
$\Def_s$ (or the formal deformation space)
of holomorphically symplectic deformations
on $M$ is equipped with the period map
$\Def_s\arrow H^2(M,\C)$, mapping
$(M,\Omega)$ to the cohomology class
of the holomorphically symplectic form
$\Omega$ in $H^2(M, \C)$.  In \cite{_KV:deformations_}
it was shown that when $H^{>0}(\calo_M)$ vanishes,
the period map is an isomorphism in the formal
category. This result holds generally,
the K\"ahler condition and compactness are not 
required. 

In \cite{_KV:deformations_},
the holomorphically symplectic deformation
space was studied using the cohomology
of the sheaf $\Ham(M)$ of holomorphic
Hamiltonian vector fields on $M$.
The corresponding DG-Lie algebra $(\Lambda^{0,*}(M)\otimes \Ham(M), \bar\6)$
is responsible for the holomorphically symplectic deformations of $(M,\Omega)$.
However, the sheaf $\Ham(M)$ can be written explicitly by
\[
0\arrow \C_M\arrow \calo_M \stackrel \Theta \arrow \Ham(M) \arrow 0,
\]
where $\C_M$ is the constant sheaf, and $\Theta(f)=\Omega^{-1}(df)$
is the Hamiltonian vector field associated with the holomorphic
function $f$. When $H^{>0}(\calo_M)=0$, the multiplication on
$H^{>0}(\Ham(M))$ vanishes because $H^i(\Ham(M))=H^{i+1}(\C_M)$
and the product $H^i(\Ham(M))\times H^j(\Ham(M))\arrow H^{i+j}(\Ham(M))$
can be represented as
\begin{equation}\label{_Ham_commu_Intro_Equation_}
H^i(\Ham(M))\times H^{j+1}(\C_M)\arrow H^{i+j+1}(\C_M),
\end{equation}
with the multiplication given by the Poisson product.
The latter clearly vanishes, because the Poisson product with the
constant vanishes. 

Pushing this argument further, in \cite{_KV:deformations_}
vanishing of obstructions of deformations was proven 
whenever the natural map $H^i(M, \C) \arrow H^i(\calo_M)$
is surjective for all $i>0$. However, even this condition
is not strong enough for our present purposes.
In this paper, we prove a version of the main result
of \cite{_KV:deformations_} which can be applied to
Bogomolov-Guan manifolds.

\hfill

\theorem (\ref{_TT_for_HS_Theorem_})\\
Let $(M,I, \Omega)$ be a compact holomorphically symplectic
manifold (not necessarily K\"ahler). Assume that the Dolbeault
cohomology group $H^{0,2}_{\bar \6}(M)= H^2(\calo_M)$ is generated
by $\6$-closed $(0,2)$-forms, and all 
$\6$-exact holomorphic 3-forms on $M$ vanish.  Then 
\begin{description}
\item[(a)]
the holomorphic symplectic
deformations of $(M,I, \Omega)$ are unobstructed.
\item[(b)]
If, in addition, all classes in the Dolbeault
cohomology group $H^{1,1}_{\bar \6}(M) $ are represented by closed $(1,1)$-forms,
then the complex deformations of $M$ are also unobstructed,
and all sufficiently small complex deformations remain
holomorphically symplectic. \endproof 
\end{description}

We prove that the Bogomolov-Guan manifolds
satisfy both assumptions of this theorem in 
\ref{_BG-conditions-BB-existence_Proposition_}.

\subsection{Beauville-Bogomolov-Fujiki form}

The Beauville-Bogomolov-Fujiki (BBF) form $q$ is a 
bilinear symmetric form on the second cohomology
of a compact hyperk\"ahler manifold which satisfies
$b_1(M)=0$ and $h^{2,0}(M)=1$. This form has topological
nature and can be defined using the Fujiki formula, 
see \eqref{_Fujiki_w^2n_Equation_} below.

Existence of the BBF form
was originally established using the deformation theory.
Bogomolov, working with compact K\"ahler holomorphically
symplectic manifolds, has noticed that the period space
of holomorphic symplectic forms $\Perspace\subset H^2(M,\C)$
has codimension 1 in $H^2(M,\C)$. Therefore, $\Perspace$
is locally identified with the set of all 
non-zero $\Omega\in H^2(M,\C)$ such that 
$P(\Omega)=0$, where 
\[
P(\Omega):\;= \int_M \Omega^{2n}=0, 
\]
and $2n=\dim_\C M$. When $H^2(M)$ admits the
Hodge decomposition, there is a $U(1)$-action $\rho$
on $H^2(M)$ and on the subalgebra on cohomology
generated by $H^2(M)$, with $\rho(t)\restrict{H^{p,q}(M)}=t^{\1(p-q)}\Id$.
Since the polynomial $P(\Omega)$ is $\rho$-invariant,
for any 2-plane $W\subset H^2(M,\R)$ where $\rho$ acts by rotations,
$P(\eta)$ is an $n$-th power of quadratic form
$q(\eta,\eta)$ on $H^2(M, \R)$.
Using the deformation theory, we produce sufficiently many
such 2-planes, proving the famous Fujiki formula:
\begin{equation}\label{_Fujiki_w^2n_Equation_}
\int_M \eta^{2n} = q(\eta,\eta)^n 
\end{equation}
where $q$ is a quadratic form on $H^2(M)$, called
{\bf Beauville-Bogomolov-Fujiki form}.
We define the BBF form for Bogomolov-Guan manifolds in
\ref{_BB_main_Theorem_}.

As a formal consequence of \eqref{_Fujiki_w^2n_Equation_},
we obtain another formula, also known as "Fujiki Formula".
It is obtained by the polarization of the $2n$-forms on both sides of
the equation \eqref{_Fujiki_w^2n_Equation_}:
\[
\int\eta_1\wedge \eta_2 \wedge ... \wedge \eta_{2n} =
2^{2n-2} (C^{2n}_n)^{-1}
\sum_{\sigma\in \Sigma_{2n}} q(\eta_{\sigma_1}, \eta_{\sigma_2})
q(\eta_{\sigma_3}, \eta_{\sigma_4})...q(\eta_{\sigma_{2n-1}}, \eta_{\sigma_{2n}})
\]
(see \cite{_Fujiki:HK_}). 
Here the sum is taken over all permutations $\sigma\in \Sigma_{2n}$.

In the present paper, we apply the deformation theory in
the case of simply-connected holomorphically symplectic
non-K\"ahler manifolds, in particular, Bogomolov-Guan
manifolds, obtaining the local Torelli theorem. This
theorem is used to deduce a version of Fujiki formula.

\hfill

\theorem (\ref{_BB_main_Theorem_})
Let $M$ be a compact holomorphically symplectic manifold,
$\dim_\C M=2n$, admitting the Hodge decomposition 
on $H^2(M)$.\footnote{\ref{_Hodge_deco_Definition_}.}
Assume that all $\6$-exact holomorphic 3-forms on $M$ vanish,
and $H^{0,2}(M)=H^{2,0}(M)=\C$. Then the space
$H^2(M)$ is equipped with a bilinear symmetric form
$q$ such that for any $\eta\in H^2(M)$, one has
$\int_M \eta^{2n} = \lambda q(\eta, \eta)^n$, where
$\lambda$ is a fixed constant. 

As a corollary of this theorem we have the BBF-form on Bogomolov-Guan manifolds (see \ref{BBF-form-BG}).


\section{Preliminaries}


\subsection{Holomorphically symplectic manifolds}

The main object of this paper is holomorphically
symplectic non-K\"ahler manifolds, which are related to
hyperk\"ahler manifolds defined below.

\hfill

\definition A holomorphically symplectic manifold 
is a complex manifold equipped with a non-degenerate, holomorphic
$(2,0)$-form $\Omega$.

\hfill

\definition (E. Calabi, \cite{_Calabi_})\\ Let $(M, g)$ be a Riemannian
manifold equipped with three complex structure operators
$I, J, K:\; TM\arrow TM$, satisfying the quaternionic relation
\[ I^2=J^2=K^2=IJK=-\Id.\]  Suppose that $I$, $J$, $K$ are
K\"ahler. Then $(M, I, J, K, g)$ is called \textbf{hyperk\"ahler manifold}.

\hfill

\remark 
The term ``hyperk\"ahler'' is often
used as a synonym of ``holomorphically symplectic''.
For compact K\"ahler manifolds
there is essentially no difference between hyperk\"ahler
geometry and holomorphic symplectic geometry. Indeed,
hyperk\"ahler manifolds are holomorphically symplectic
since $\Omega:\;=\omega_J+\1\omega_K$ is a holomorphic
symplectic form on $(M,I)$. Converse is also true in the K\"ahler case by
Calabi-Yau theorem:

\hfill

\theorem (Calabi-Yau; \cite{_Yau:Calabi-Yau_})\\
Let $M$ be a compact, holomorphically symplectic K\"ahler
manifold. Then $M$ admits a hyperk\"ahler metric, which is
uniquely determined by the cohomology class of its 
K\"ahler form $\omega_I$.

\hfill

\definition A hyperk\"ahler manifold $M$ is called
{\bf of maximal holonomy}, or {\bf simple},
or {\bf IHS}, if $\pi_1(M)=0$, $H^{2,0}(M)=\C$.

\hfill

Simple hyperk\"ahler manifolds are building blocks of all compact hyperk\"ahler manifolds due to the famous Bogomolov's theorem.
 
\hfill

\theorem (Bogomolov's decomposition: \cite{_Bogomolov:decompo_})\\ 
Any hyperk\"ahler manifold admits a finite covering
which is a product of a torus and several 
maximal holonomy hyperk\"ahler manifolds.
The maximal holonomy hyperk\"ahler components
of this decomposition are defined uniquely.

\hfill

There are several examples of hyperk\"ahler manifolds are known: two infinite series -- Hilbert scheme $K3^{[n]}$ of K3, and generalized Kummer manifolds, and beside of this, two sporadic examples -- six-, and ten-dimensional manifolds of O'Grady type.

\subsection{Deformation theory of holomorphically symplectic manifolds} 
\label{def-theory-subsection}

Now we will briefly introduce the
main notions of the deformation theory, for details see \cite{_Angela_}.

\hfill

\proposition
Every compact complex manifold $X$ admits a deformation $\pi:\; X_U \to \Def(X)$, where $\Def(X)$ is a germ of an analytic set, such that for any deformation $\pi_0:\; X \to S$ there exists a morphism $g:\; (\Def(X), 0) \to (S, 0)$ such that $X$ is isomorphic to the pullback $X_U \times_{\Def(X)} S$
and the differential $dg_0:\; T_0 \Def(X) \to T_0 S$ is uniquely determined. 

\hfill

\definition This deformation is called the {\bf Kuranishi deformation} of $X$. If
moreover for any $\pi_0$ such a morphism $g$ is unique,
then the deformation is called {\bf universal}.

\hfill

The original Bogomolov-Tian-Todorov theorem says that deformation theory works 
when $M$ has trivial canonical bundle and $\6 \bar \6$-lemma (see \cite{_Bogomolov:TT_,_Tian:TT_, _Todorov:TT_}) for original
proofs.

\hfill

\theorem (Bogomolov-Tian-Todorov theorem) Let $X$ be a compact
complex manifold such that its canonical bundle is holomorphically trivial. Suppose that $X$ satisfies the $\6 \bar \6$-lemma.
Then $X$ has unobstructed deformations and
the base of its Kuranishi deformation is smooth.

\hfill

The proof of this theorem uses the Tian-Todorov lemma 
(\ref{_TT_Lemma_}), which holds true on all Calabi-Yau
manifolds, and $\6\bar\6$-lemma, which is true on
compact K\"ahler manifolds.  
Iacono \cite{_Iacono_} showed that the $\6\bar\6$-lemma
can be replaced by a more abstract and general condition,
the ``homotopy abelianness'' of the differential graded Lie
algebra (DGLA) $(\Lambda^{0,*}(TM), \bar\6)$ which is responsible
for the deformations. 

For Bogomolov-Guan manifolds, 
the $\6\bar\6$-lemma is most likely false.
Indeed, it is known to be false
for the Kodaira-Thurston surfaces, since $b_1 \neq 2 h^{0,1}$ (\cite{G}). It seems
very unlikely that the traditional approach
to deformations of Calabi-Yau would work
for the Bogomolov-Guan manifolds.

Kaledin and Verbitsky in \cite{_KV:deformations_}
suggested an alternative approach to deformations
of holomorphically symplectic manifolds related to
\cite{_Barannikov_Kontsevich_},
\cite{_Kontsevich:lectures_}. They study
holomorphically symplectic deformations, i.e. deformations of the holomorphically symplectic
form (the complex structure is determined
by the holomorphically symplectic form uniquely).
In this case the Kuranishi deformation space
$\Def_s$
of holomorphically symplectic deformations
on $M$ is equipped with the period map
$\Def_s\arrow H^2(M,\C)$, mapping
$(M,\Omega)$ to the cohomology class
of the holomorphically symplectic form
$\Omega$ in $H^2(M, \C)$.  It was shown that when $H^{>0}(\calo_M)$ vanishes,
the period map is an isomorphism in the formal
category. This result (\cite[Theorem 1.1]{_KV:deformations_}) does not require 
the K\"ahler condition and compactness.

In the present paper, we modify the Kaledin-Verbitsky
method to prove the local Torelli theorem for
the Bogomolov-Guan manifolds.

\subsection{Beauville-Bogomolov-Fujiki form} 
\label{BBF-subsection}

Let $M$ be a hyperk\"ahler manifold of maximal holonomy
(that is, satisfying $\pi_1(M)=0$ and $H^{2,0}(M)=\C$).
Then the space $H^2(M,\Q)$ is equipped with a non-degenerate bilinear
symmetric quadratic form $q$ which is called
Beauville-Bogomolov-Fujiki form (BBF-form). Existence of
this form in the K\"ahler case can be deduced from the
Fujiki theorem.

\hfill

\theorem
(Fujiki, \cite{_Fujiki:HK_})\\ 
Let $\eta\in H^2(M)$, and $\dim M=2n$, where $M$ is
hyperk\"ahler. Then there is primitive integer quadratic form $q$ on $H^2(M,\Z)$, and $c>0$ an integer number, such that for any two-form $\eta$
the following holds
$\int_M \eta^{2n} = q(\eta,\eta)^n$

\hfill

\definition \label{_BBF_Definition_}
This form is called
{\bf the Beauville-Bogomolov-Fujiki form}.  It is defined
by the Fujiki's relation uniquely up to a sign. The sign is determined
from the following formula (Bogomolov, Beauville)
\begin{align*}  \lambda q(\eta,\eta) &=
   \int_X \eta\wedge\eta  \wedge \Omega^{n-1}
   \wedge \bar \Omega^{n-1} -\\
 &-\frac {n-1}{2n}\left(\int_X \eta \wedge \Omega^{n-1}\wedge \bar
   \Omega^{n}\right) \left(\int_X \eta \wedge \Omega^{n}\wedge \bar \Omega^{n-1}\right)
\end{align*}
where $\Omega$ is the holomorphic symplectic form, and 
$\lambda>0$.

\hfill

\remark The form $q$ has signature $(b_2-3,3)$.
It is negative definite on primitive forms, and positive
definite on $\langle \Omega, \bar \Omega, \omega\rangle$,
 where $\omega$ is a K\"ahler form.

\hfill

\remark The local Torelli theorem and Fujiki formula
also hold for some singular holomorphic
symplectic varieties (\cite{_Bakker_Lehn_,_Kirshner,_Namikawa_, _Namikawa_form}).

\subsection{Bogomolov-Guan manifolds} 
\label{BG-manifolds-subsection}

In this section we describe the construction of
Bogomolov-Guan manifolds, the
non-K\"ahler simply-connected holomorphically symplectic
manifolds which are the cental object of this paper.


\subsubsection{Kodaira-Thurston surfaces}

\definition \cite[p. 197]{_Barth_Peters_Van_de_Ven_} \label{KT_surf}
Let $L$ be a line bundle
on an elliptic curve $E$ with the first Chern class
$c_1(L)\neq 0$. Denote by $\tilde S$ the corresponding $\C^*$-bundle
on $E$ obtained by removing the zero section,
$\tilde S=\Tot(L) \backslash 0$. Fix a complex
number $\lambda$ with $|\lambda| > 1$, and
let $h_\lambda:\; \tilde S \arrow \tilde S$
be the corresponding homothety of $\tilde S$.
The quotient $\tilde S/\langle h_\lambda\rangle$
is called {\bf a (primary) Kodaira-Thurston surface}.
This surface was well known in algebraic geometry since
the mid-1960ies. 
It first appeared in the paper \cite{_Kodaira:surfaces_1_}
as Kodaira classified complex surfaces with trivial
canonical bundle. Kodaira proved that such surfaces
are isomorphic to a torus, K3 surface or a Kodaira
surfaces, and proved that the Kodaira surface
has algebraic dimension 1 (\cite[Theorem 19]{_Kodaira:surfaces_1_}).

In his 1976 paper \cite{_Thurston:Kodaira_}, 
Thurston independently constructed this
manifold and proved that it is symplectic but does not
support a K\"ahler structure. The name
``Kodaira-Thurston'' originates in mid-1980-es due to the 
Spanish mathematicians who studied the topology and
differential geometry of the complex and
symplectic nilmanifolds (\cite{_de_Leon:onThurston_}, 
\cite{_AFM:KT_}).

\hfill

\remark \label{KT_nilman}
We could also consider Kodaira-Thurston surface as
follows: $S_r = G/Z_r,$
where
\begin{equation}\label{_KT_nilma_Equation_}
G = (x,y) \in  \left( \begin{array}{ccc}
1 & x & y \\
0 & 1 & \bar x \\
0 & 0 & 1 \end{array} \right), \quad Z = \lbrace (x,y) \in
G|_{x\in r(\Z+i\Z), y \in r(\Z+i\Z) },  r \in \Z \rbrace 
\end{equation}
This is a holomorphically symplectic manifold with a 2-form $\omega = dx \wedge dy$.

\hfill

Guan's construction of higher-dimensional analogues \cite{Gu2} uses
this definition of a Kodaira-Thurston surface as a starting point.

\hfill

\claim \label{KT_proj_to_torus}
The nilmanifold \eqref{_KT_nilma_Equation_} is isomorphic to the Kodaira
surface defined in  \ref{KT_surf}. 

\hfill

\proof
Consider the line bundle $L$ with $c_1(L) \neq 0$, 
then $S$ is a compact surface which is fibered with
elliptic fiber $\mathbb{C}^*/Z = E_L$ over
$E$, where $Z \subset \mathbb{C}^*$  is
a discrete cocompact subgroup. The canonical
bundle of $S$ is trivial, and the global holomorphic 
form $\Omega$ on $S$ is invariant under the action of $E_L$.
Topologically, $S$ has a structure of a principal
fibration $p_T$ over $T^3$ with $S^1$ as a fiber. A
projection $p_s:\; T^3=E\times S_0^1 \rightarrow E$
contracting $S_0^1$ satisfies  $p_E = p_s p_T$, where
$p_E$ is a principal holomorphic $E_L$ fibration over $E$.
\endproof

\hfill

These surfaces have $b_1 = 3$, therefore they are not K\"ahler.

\subsubsection{Construction of Bogomolov-Guan manifolds} \label{BG-construction_section}

Higher-dimensional simply-connected non-K\"ahler
holomorphically symplectic manifolds have been discovered by
Guan (\cite{Gu1, Gu2, Gu3}) and by Bogomolov
(\cite{B1}). Below we briefly describe Bogomolov's construction.

Construction begins with a surface $S$ with $c_1 ( L ) = m
> 2$. Let $n > 2$ be a
number dividing $m$. From this data we construct a
holomorphically symplectic manifold
of dimension $2n - 2$. We follow the Bogomolov's approach.

Consider the $n$-symmetric power $S^{(n)}$ of $S$. It is
a singular complex variety admitting the standard Douady
desingularization, $S^{[n]}$, called the Hilbert scheme
of points on $S$. It is well-known (see
\cite{_Beauville_}) that $S^{[n]}$ is holomorphically
symplectic whenever $S$ is holomorphically
symplectic. There exists a holomorphic projection $\pi_e:\; S^{(n)}
\rightarrow E$ induced from the elliptic fibration
$\pi_E:\; S \arrow E$ such that and $\pi_e=p_s\pi_t$,
where $\pi_t:\; S^{(n)} \rightarrow T^3$ is the map induced by $p_T$ (see
\ref{KT_proj_to_torus}).

Taking the composition with 
the resolution $r:\; S^{[n]} \arrow S^{(n)}$,
we obtain an isotrivial fibration
$\pi:\; S^{[n]}\arrow E$ and 
$\sigma:\; S^{[n]}\arrow T^3$ satisfying $\pi = p_s \sigma$. Denote a fiber of $\pi$ by $F^{[n]}$.
Then $F^{[n]}$ is a smooth divisor in 
a holomorphically symplectic manifold $(S^{[n]}, \Omega)$ with locally free action of $E_L$.
The restriction of $\Omega$ to $F^{[n]}$ has rank 
$2n-2$, because $F^{[n]}\subset S^{[n]}$ is a divisor.

\hfill

\definition \label{def_charac_foliation} Let $D$ be a divisor in
holomorphically symplectic manifold. Denote by $K \subset TD$ the kernel of $\Omega\restrict {D}$,
that is, the set of all $x\in TD$ such that 
$\Omega\restrict {D}(x, \cdot)=0$.
This foliation is called {\bf the characteristic foliation}.

\hfill

\claim
Let $S$ be a Kodaira surface, $\pi:\; S^{[n]}\arrow E$ the isotrivial
fibration constructed above, and $F^{[n]}$ its fiber. Then
the leaves of  its characteristic foliation $K$ coincide with the orbits 
of $E_L$-action on $F^{[n]}$. 

\proof Consider the fibration $\tilde \pi:\; S^n\arrow
E$. Clearly, the leaves of the characteristic foliation on 
$\tilde\pi^{-1}(0)$ coincide with the $E_L$-orbits. Clealry,
$S^{[n]}$ is obtained from $S^n$ by taking a finite
quotient and desingularizing. Both of these operations
are compatible with taking the characteristic
foliation. \endproof

\hfill

\theorem \cite[Lemma 3.13]{B1} \label{BG_existence} There
exists a compact complex simply-connected holomorphically
symplectic manifold $Q$ of $\dim_\C=2n-2$ with the
surjective map $p_w :\; Q \rightarrow W$ of degree 
$n^2$ that $p_w^* \omega_0$ is a non-degenerate closed holomorphic (2,0)-form
on $Q$.

\hfill

\definition \label{_Bogomolov_Guan_Definition_}
The smooth order $n^2$ covering $Q$  of 
$W$ is called \textbf{the Bogomolov-Guan manifold}.

\hfill

We give an explicit construction of $Q$  as
follows. Denote $\sigma^{-1}(0)$ by $H^{[n]}$,
where $\sigma:\; S^{[n]}\arrow T^3$ is the 
projection to the 3-dimensional torus $T^3$ defined above.
Clearly, $\sigma^{-1}(0)$ is a
real manifold of dimension $4n-3$. Consider the quotient
$R := H^{[n]} / S^1$. It is a compact complex
variety with quotient singularities and a map 
$p :\; R \to W$ of degree $n$ (\cite{B1}; see the the diagram
\eqref{BG_construction} below). The
manifold $Q$ can be obtained as a cyclic $n$-covering $p_r$ of $R$
\cite[Lemma 3.6, Corollary 3.7]{B1}; it is the universal
cover of $R$ considered as an orbifold. 

The following diagram relates the varieties used in this construction.
\begin{equation} \label{BG_construction}
\begin{diagram}[labelstyle=\tiny]
T & & & &H^{[n]}&\rTo^{S^1}&R &&\\
& \luTo^{\sigma}& & \ldInto & & & &\luTo^{p_r} &\\
\dTo_{p_s}^{S^1} & & S^{[n]}& & && \dTo_{p} & & Q\\
& \ldTo_{\pi} & & \luInto & & & &\ldTo_{p_w}&\\
E & && &F^{[n]}& \rTo_{E_L} & W&&\\ \end{diagram}
\end{equation}


The non-K\"ahlerness follows from the existence of a nonsingular
surface $V \subset Q$ isomorphic to 
a blow-up of the Kodaira surface $S$ \cite[Corollary 4.10]{B1}.

\section{Deformations of holomorphically symplectic manifolds}


\subsection{Tian-Todorov lemma for holomorphically symplectic manifolds}

\definition
Let $M$ be a complex manifold, and $\Lambda^{0,p}(M)\otimes T^{1,0}M$ the sheaf of
$T^{1,0}M$-valued $(0,p)$-forms. Consider the commutator bracket $[\cdot, \cdot]$ on $T^{1,0}M$,
and let $\bar \calo_M$ denote the sheaf of antiholomorphic functions.
Since  $[\cdot, \cdot]$ is $\bar \calo_M$-linear, it is naturally extended
to 
\begin{equation}\label{_Schouten_Equation_}
\Lambda^{0,p}(M)\otimes_{C^\infty M} T^{1,0}M=\overline
       {\Omega^pM}\otimes_{\bar\calo_M}T^{1,0}M,
\end{equation}
giving a bracket 
\[ [\cdot, \cdot]:\; \Lambda^{0,p}(M)\otimes T^{1,0}M\times \Lambda^{0,q}(M)\otimes T^{1,0}M
\arrow \Lambda^{0,p+q}(M)\otimes T^{1,0}M.
\]
This bracket is called {\bf the Schouten bracket}. 

\hfill

\remark
Since $[\cdot, \cdot]$ is $\bar \calo_M$-linear,
the Schouten bracket satisfies the Leibnitz identity:
\[
\bar\6([\alpha, \beta])= [\bar\6\alpha, \beta]+[\alpha, \bar\6\beta].
\]
This allows one to extend the Schouten bracket to the
$\bar\6$-cohomology of the complex
$(\Lambda^{0,*}(M)\otimes T^{1,0}M, \bar\6)$,
which coincide with the cohomology of the sheaf of holomorphic vector fields:
$ [\cdot, \cdot]:\; H^p(TM) \times H^q(TM) \arrow H^{p+q}(TM)$.

\hfill

Let now  $\Omega$ be a holomorphically symplectic form on a complex manifold $M$,
$\dim_\C M=2n$. Then $TM\cong \Omega^1M$, hence the Schouten bracket is defined on
the Dolbeault cohomology: 
\begin{equation}\label{_Schouten_holo_symp_Equation_}
[\cdot, \cdot]:\; H^{p,1}(M) \times H^{q,1}(M) \arrow H^{p+q,1}(M).
\end{equation}
Tian-Todorov lemma expresses the Schouten bracket on a Calabi-Yau manifold
in terms of $\6, \bar\6$-differentials on the $(p, q)$-forms.
Here we give the usual Tian-Todorov lemma and then 
a holomorphically symplectic version of 
Tian-Todorov lemma, obtaining an even simpler expression.

\hfill

Let us recall the usual Tian-Todorov lemma.
Assume that $M$ is a complex $n$-manifold with trivial canonical bundle $K_M$,
and $\Phi$ a non-degenerate section of $K_M$.
We call a pair $(M, \Phi)$ {\bf a Calabi-Yau manifold}.
Substitution of a vector field into $\Phi$
gives an isomorphism $TM \cong \Omega^{n-1}(M)$.
Similarly, one obtains an isomorphism
\begin{equation}\label{_polyvector_forms_Equation_}
 \Lambda^{0,q}M \otimes \Lambda^p TM \arrow \Lambda^{0,q}M \otimes \Lambda^{n-p,0}M =\Lambda^{n-q,p}M.
\end{equation}
{\bf Yukawa product} 
$\bullet:\; \Lambda^{p,q}M \otimes \Lambda^{p_1,q_1}M
   \arrow  \Lambda^{p+p_1 -n,q+q_1}M$ 
is obtained from the usual product
\[
 \Lambda^{0,q}M \otimes \Lambda^p TM \times \Lambda^{0,q_1}M \otimes \Lambda^{p_1} TM\arrow
\Lambda^{0,q+q_1}M \otimes \Lambda^{p+p_1} TM 
\]
using the isomorphism \eqref{_polyvector_forms_Equation_}.

\hfill

\lemma  \label{_TT_Lemma_}
(Tian-Todorov lemma)\\
Let $(M,\Phi)$ be a Calabi-Yau manifold, and 
\[ [\cdot, \cdot]:\; \Lambda^{0,p}(M)\otimes T^{1,0}M\times \Lambda^{0,q}(M)\otimes T^{1,0}M
\arrow \Lambda^{0,p+q}(M)\otimes T^{1,0}M.
\]
its Schouten bracket. Using the isomorphism
\eqref{_polyvector_forms_Equation_}, we can interpret
Schouten bracket as a map
\[ [\cdot, \cdot]:\; \Lambda^{n-1,p}(M)\times \Lambda^{n-1,q}(M)
\arrow \Lambda^{n-1,p+q}(M).
\]
Then, for any $\alpha \in \Lambda^{n-1,p}(M)$, $\beta \in \Lambda^{n-1,p_1}(M)$, one has
\begin{equation}\label{_TT_Equation_}
[\alpha, \beta] = \6(\alpha \bullet \beta) - (\6\alpha)\bullet \beta - (-1)^{n-1+p}\alpha \bullet (\6\beta), 
\end{equation}
where $\bullet$ denotes the Yukawa product.

\hfill

\proof 
See \cite{_Tian:TT_}, \cite{_Todorov:TT_}. 
\endproof

\hfill

\definition
Let now $(M, \Omega)$ be a holomorphic symplectic manifold again, $\dim_\C M =2n$.
Holomorphic symplectic form gives a pairing
$\Omega^1 M \otimes \Omega^1 M \arrow \calo_M$.
Extend this pairing to 
$\Omega^i M \otimes \Omega^i M \arrow \calo_M$
by multilinearity, $\alpha, \beta \arrow (\alpha, \beta)_\Omega$. Define the holomorphic
symplectic $\star$-map 
\[ \star: \; \Lambda^{p,0}(M)\arrow \Lambda^{n-p,0}(M)
\]
via 
\[
(\alpha, \beta)_\Omega= \frac{\alpha \wedge \star \beta}{\Omega^n}.
\]
This is the usual Hodge star operator on $(1,0)$-variables,
with the holomorphic volume form used instead of the usual volume form.
We extend $\star$-map to $\Lambda^{p,q}(M)$ by
$\star(\alpha\wedge\gamma)= \star(\alpha)\wedge\gamma$
for any $(0,p)$-form $\gamma$.

\hfill

\lemma
Let $M$ be a holomorphic symplectic manifold.
Consider the operators $L_\Omega(\alpha):= \Omega\wedge \alpha$,
$H_\Omega$ acting as multiplication by $n-p$ on $\Lambda^{p,q}(M)$,
and $\Lambda_\Omega:= \star \Lambda\star$.
Then $L_\Omega, H_\Omega, \Lambda_\Omega$ satisfy the
$\goth{sl}(2)$ relations, similar to the Lefschetz triple:
\[
[H_\Omega, L_\Omega]=2 L_\Omega, \ \ \ [H_\Omega, \Lambda_\Omega]=-2 \Lambda_\Omega, 
[L_\Omega, \Lambda_\Omega]=H_\Omega.
\]
\proof See \cite[Section 9]{_V:Mirror_}.
\endproof

\hfill

Now we can state the holomorphically symplectic version of Tian-Todorov lemma.

\hfill

\theorem\label{_hol_symp_TT_Theorem_}
Let $(M, \Omega)$ be a holomorphically symplectic manifold, and
\[ [\cdot, \cdot]:\; \Lambda^{0,p}(M)\otimes T^{1,0}M\times \Lambda^{0,q}(M)\otimes T^{1,0}M
\arrow \Lambda^{0,p+q}(M)\otimes T^{1,0}M.
\]
the Schouten bracket. Using $\Omega$, we identify $T^{1,0}(M)$ and $\Lambda^{1,0}(M)$,
and consider the Schouten bracket as a map
\begin{equation}\label{_Sch_in_hol_symp_Equation_}
[\cdot, \cdot]_\Omega:\; \Lambda^{1,p}(M)\times \Lambda^{1,q}(M)
\arrow \Lambda^{1,p+q}(M).
\end{equation}
Then for any $a, b \in \Lambda^{1,*}(M)$, one has
\[
[a, b] = \delta(a \wedge b) - (\delta a) \wedge b - (-1)^{\tilde a}a\wedge\delta(b),
\]
where $\tilde a$ is parity of $a$, and $\delta:= [\Lambda_\Omega, \6]$.

\hfill

\proof
Acting by $\star$ on $\Lambda^{*,*}(M)$, we obtain the 
Yukawa multiplication from the usual multiplication:
$\alpha \bullet\beta=\pm \star(\star \alpha\wedge \star \beta)$. Moreover,
the Schouten bracket interpreted as in Tian-Todorov lemma gives a map
\[ [\cdot, \cdot]:\; \Lambda^{n-1,p}(M)\times \Lambda^{n-1,q}(M)
\arrow \Lambda^{n-1,p+q}(M).
\]
after twisting by $\star$ becomes the bracket $[\cdot,
  \cdot]_\Omega$
defined in \eqref{_Sch_in_hol_symp_Equation_}.
Then \eqref{_TT_Equation_} becomes
\[
[a, b]_\Omega = 
\star \6 \star (a \wedge b) - (\star \6 \star a) \wedge b - 
(-1)^{\tilde a}a\wedge\star \6 \star(b).
\]
Therefore, \ref{_hol_symp_TT_Theorem_} would be implied if we prove
the following holomorphic symplectic analogue of the K\"ahler relations:
\begin{equation} \label{_star_6_via_commu_Equation_}
\star \6 \star=[\Lambda_\Omega, \6].
\end{equation} 

This statement is local, hence it would suffice to prove it in a coordinate patch.
Choose holomorphic Darboux coordinates such that 
$\Omega=\sum_{i=1}^{n} dz_{2i-1}\wedge dz_{2i}$, let $\alpha$ 
be a coordinate monomial and $f$ a function.
A simple computation gives 
$\delta(f\alpha)= \sum_i \frac{\6 f}{\6 z_i} i_{\Omega^{-1}(dz_i)}(\alpha)$,
where $i_v(\alpha)$ is convolution of a vector field $v$ and a form $\alpha$, and
$\Omega^{-1}(dz_i)$ the vector field dual to $dz_i$ via $\Omega$.
Similarly, $\star \6 \star (f\alpha)= \sum_i \frac{\6 f}{\6 z_i} \star e_{dz_i}\star \alpha$,
where $e_{dz_i}(\alpha)=dz_i \wedge \alpha$. Therefore,
\eqref{_star_6_via_commu_Equation_} and \ref{_hol_symp_TT_Theorem_}
follow from $\star e_{dz_i} \star = i_{\Omega^{-1}(dz_i)}$.
This is clear, because $\star e_{dz_i} \star$ is the convolution
with a vector field $\Omega$-dual to $dz_i$.
\endproof

\subsection{Deformation of complex manifolds and the Schouten bracket}
\label{_obstructions_def_Subsection_}

Here we reproduce the standard construction of the deformation
theory for complex structures; we follow 
\cite{_Barannikov_Kontsevich_}, \cite{_Kontsevich:lectures_}
and \cite{_Todorov:TT_}. We deduce the vanishing of obstructions
from the Tian-Todorov lemma, proving the Bogomolov-Tian-Todorov
theorem on unobstructedness of the deformations of Calabi-Yau manifolds.

Let $(M,I)$ be an almost complex manifold, and $B$ an abstract vector bundle
over $\C$ isomorphic to $\Lambda^{0,1}(M)$. To identify $B$ and $\Lambda^{0,1}(M)$
we can use a first order differential operator $\bar\6:\; C^\infty M \arrow B=\Lambda^{0,1}(M)$
satisfying the Leibnitz rule. Using the universal property of K\"ahler differentials,
we obtain that any such operator is associated with a linear map
$u:\; \Lambda^1(M,\C)\arrow B$. Then $B= \frac{\Lambda^1(M,\C)}{\ker u}= \Lambda^{0,1}(M)$.
The almost complex structure operator on $M$ is reconstructed 
as follows: we represent $\Lambda^1_\C M$ as a direct sum 
$\Lambda^1_\C M=\ker u \oplus \overline{\ker u}$, and write $I$
as $\1$ on $\ker u$ and $-\1$ on $\overline{\ker u}$.
The main advantage of this approach is that the
integrability of almost complex structure is very easy to write.
Indeed, let us extend $\bar\6:\; C^\infty M \arrow B$ to the
corresponding exterior algebra using the Leibnitz rule:
\[
C^\infty M \stackrel {\bar\6}\arrow B \stackrel {\bar\6}\arrow\Lambda^2 B
\stackrel {\bar\6}\arrow\Lambda^3 B \stackrel {\bar\6}\arrow...
\]
Then integrability of $I$ is equivalent to $\bar\6^2=0$. 
This allows one to describe the deformations of the complex
structure as follows. Let $\gamma\in TM \otimes B$ be a $B$-valued derivation
of $C^\infty M$. Then $\bar\6+ \gamma$ also satisfies the
Leibnitz rule. However, $\ker u$ does not change if
$\gamma \in T^{1,0}M \otimes B$, hence we may assume that
$\gamma\in T^{1,0}M \otimes B= T^{1,0}M \otimes \Lambda^{0,1}(M)$.
Therefore, almost complex deformations of $I$ are given by the
sections $\gamma\in T^{1,0}M \otimes \Lambda^{0,1}(M)$,
with the integrability relation $(\bar \6+\gamma)^2=0$
rewritten as Maurer-Cartan equation\footnote{Depending on
  conventions, the Maurer-Cartan equation can be written
as $\bar \6(\gamma) = -\frac{1}{2} \{\gamma, \gamma\}$. This version
is obtained from the one we use by rescaling.}
\begin{equation}\label{_Maurer_Cartan_Equation_}
\bar \6(\gamma) = -\{\gamma, \gamma\}.
\end{equation}
Here $\bar \6(\gamma)$ is identified with the anticommutator
$\{\bar \6, \gamma\}$, and $\{\gamma, \gamma\}$ is anticommutator
of $\gamma$ with itself, where $\gamma$ is considered as a $\Lambda^{0,1}(M)$-valued
differential operator. This identifies $\{\gamma, \gamma\}$ with the Schouten bracket.
Further on, we shall write $[\gamma, \gamma]$ instead of $\{\gamma, \gamma\}$, 
because this usage is more common.

The Kuranishi deformation space of complex structures on $M$
is identified with the space of solutions of \eqref{_Maurer_Cartan_Equation_}
modulo the diffeomorphism action. Here we are interested in non-obstructedness
of these deformations, which is defined as follows.

Let us write $\gamma$ as power series,
$\gamma= \sum_{i=0}^\infty t^{i+1} \gamma_i$.
Then \eqref{_Maurer_Cartan_Equation_} becomes
\begin{equation}\label{_MC_series_Equation_}
\bar\6\gamma_0=0,\ \ \ \bar \6\gamma_p = -\sum_{i+j=p-1} [\gamma_i, \gamma_j].
\end{equation}
The following exact sequence of sheaves gives an acyclic resolution
for the sheaf $TM$ of holomorphic vector fields:
\[
0 \arrow TM \stackrel {\bar\6}\arrow TM \otimes \Lambda^{0,1}(M) \stackrel {\bar\6}
\arrow TM \otimes \Lambda^{0,2}(M)
\stackrel {\bar\6}\arrow TM \otimes \Lambda^{0,3}(M) \stackrel {\bar\6}\arrow....
\]
This allows one to define the cohomology class $[\gamma_0]\in H^1(M, TM)$
for any $\gamma_0$ with $\bar\6\gamma_0=0$.
We say that deformations of complex structures are
{\bf unobstructed} if the solutions $\gamma_1, ..., \gamma_n, ... $
of \eqref{_MC_series_Equation_} can be found for $\gamma_0$ in any given
cohomology class $[\gamma_0]\in H^1(M, TM)$.

Notice that the sum $\sum_{i+j=p-1} [\gamma_i, \gamma_j]$ is always $\bar\6$-closed.
Indeed, the Schouten bracket commutes with $\bar\6$, hence
\begin{equation}\label{_triple_commu_Equation_}
\bar\6\left(\sum_{i+j=p-1} [\gamma_i, \gamma_j]\right)
= \sum_{i+j+k=p-1}[\gamma_i, [\gamma_j, \gamma_k]]+ [[\gamma_i, \gamma_j], \gamma_k].
\end{equation}
vanishes as a sum of triple commutators.
Therefore, obstructions to deformations are given by
cohomology classes of the sums $\sum_{i+j=p-1} [\gamma_i, \gamma_j]$,
which are defined inductively. These classes are called
{\bf Massey powers} of $\gamma_0$. Indeed, one could define
Massey products in a differential graded Lie algebra $A$ in 
terms of the obstructions to the solution of 
\eqref{_MC_series_Equation_} in the algebra of upper
triangular matrices over $A$ (\cite{_Babenko_Taimanov_}).

In the K\"ahler setting, Tian-Todorov lemma immediately implies the
unobstructedness of deformations for compact manifolds with trivial canonical bundle. Indeed, 
we can always start from $\gamma_0 \in TM \otimes \Lambda^{0,1}(M)= \Lambda^{n-1,1}(M)$
which is harmonic. Then it is $\bar\6$- and $\6$-closed. Therefore,
$[\gamma_0, \gamma_0] = \6(\gamma_0\bullet\gamma_0)$ is $\6$-exact.
It is also $\bar\6$-closed, because the
Yukawa product commutes with $\bar\6$. Then $\6\bar\6$-lemma
implies that $[\gamma_0, \gamma_0]$ is $\6\bar\6$-exact.
Using induction, we may assume that the solutions
of \eqref{_MC_series_Equation_} for
$p=1, ..., n-1$ are all $\6\bar\6$-exact.
To solve \eqref{_MC_series_Equation_} for $p=n$, we use Tian-Todorov lemma again, 
obtaining an equation
\[
\bar \6\gamma_n = \sum_{i+j=n-1} \6(\gamma_i\bullet\gamma_j).
\]
This sum is $\bar\6$-closed by \eqref{_triple_commu_Equation_} 
and $\6$-exact, hence (by $\6\bar\6$-lemma) also $\bar\6$-exact.
This gives a solution of \eqref{_MC_series_Equation_} for any $p$.

Notice that this argument seriously depends on $\6\bar\6$-lemma.
Indeed, for non-K\"ahler manifolds with trivial canonical
bundle, deformations can be obstructed (see \cite{_Ghys_} 
for such an example).

\subsection{Hamiltonian vector fields and holomorphic symplectic deformations}
\label{_Hami_defo_Subsection_}

We apply the solution of the Maurer-Cartan equation
given in Subsection \ref{_obstructions_def_Subsection_}, 
to holomorphically symplectic manifolds. As was already noticed
in \cite{_KV:deformations_}, in this case the Tian-Todorov argument
simplifies significantly. Here we explain the relation
between the deformation theory and the cohomology of holomorphic
Hamiltonian vector fields proven in \cite{_KV:deformations_}.

Let $(M, I, \Omega)$ be a complex manifold equipped with a holomorphically
symplectic form, $B$ an abstract vector bundle
isomorphic to $\Lambda^{0,1}(M)$, 
and $\bar \6 + \gamma$ a deformation of the complex structure, 
considered in Subsection \ref{_obstructions_def_Subsection_}.
The deformation $\bar \6 + \gamma$ remains holomorphically
symplectic whenever $\Omega$ is $\gamma$-invariant, that is,
when the map $\Lambda^{2,0}(M) \stackrel\gamma \arrow \Lambda^{2,1}(M)$
defined by $\gamma$ considered as a first order differential operator
satisfies $\gamma(\Omega)=0$. This is the same as to say that
$\gamma$ is a (0,1)-form with coefficients in Hamiltonian vector fields.

Unfortunately, this is not very useful, because the
sheaf of Hamiltonian vector fields is not coherent, and one
cannot take the tensor product. 
To make this observation rigorous, we notice that the
real analytic functions are obtained by taking a completion of the tensor
product of holomorphic and antiholomorphic functions over constants.
Then $TM \otimes \Lambda^{0,1}(M)$ can be obtained by taking
a completed tensor product $\hat \otimes$ of the sheaf ${\cal T}$ of holomorphic
vector fields and the sheaf $\overline{\Omega^1M}$  of 
antiholomorphic 1-forms over the constant sheaf $\C_M$. 
This gives the definition of the relevant sheaf: 
\[
\Lambda^{0,1}(M)\otimes \Ham_M:= {\cal H}\hat \otimes_{\C_M} \overline{\Omega^1M}
\subset \Lambda^{0,1}(M)\otimes TM = {\cal T}\hat \otimes_{\C_M} \overline{\Omega^1M},
\]
where ${\cal H}$ is the sheaf of holomorphic Hamiltonian vector fields.
This is a real analytic version of the sheaf $\Lambda^{0,1}(M)\otimes \Ham_M$;
the smooth version is obtained by taking $\otimes_{\calo_\R M}C^\infty M$,
where $\calo_\R M$ is the sheaf of real analytic functions.
Abusing the language, we use the same notation for these sheaves.

A solution of the Maurer-Cartan equation
$\left(\bar \6 + \sum_{i=0}^\infty t^{i+1} \gamma_i\right)^2=0$
gives a holomorphically symplectic deformation whenever all $\gamma_i$
belong to $\Lambda^{0,1}(M)\otimes \Ham_M$.

Using $\Omega$ to dualize, the sheaf of Hamiltonian vector fields
can be embedded to $\Lambda^{1,0}(M)$ as a sheaf of $\6$-closed
(1,0)-forms.

Similarly, if we use $\Omega$ to consider $\gamma_i$ as sections of
$\Lambda^{0,1}(M) \otimes T^{1,0}M= \Lambda^{1,1}(M)$, the condition
$\gamma_i \in \Lambda^{0,1}(M)\otimes \Ham_M$ is interpreted 
as $\6\gamma_i=0$ (check \cite{_KV:deformations_} for
details). This gives the following assertion.

\hfill

\claim \label{_closed_gamma_holo_symp_Claim_}
Let $M$ be a holomorphically symplectic manifold, and
$\gamma=\sum_{i=0}^\infty t^{i+1} \gamma_i$ be a solution of Maurer-Cartan
equation $\left(\bar \6 + \sum_{i=0}^\infty t^{i+1} \gamma_i\right)^2=0$.
Using the holomorphically symplectic form as above,
we interpret $\gamma_i$ as sections of $\Lambda^{1,1}(M)$.
Suppose that  $\6\gamma_i=0$. Then the deformation
defined by $\gamma$ is holomorphically symplectic.
\endproof

\subsection{Unobstructedness of holomorphically symplectic deformations}

\definition
Let $(M,I, \Omega)$ be a holomorphically symplectic manifold.
We say that the {\bf holomorphic symplectic deformations of $(M,I, \Omega)$
are unobstructed} if for any $\bar\6$- and $\6$-closed $\gamma_0\in \Lambda^{1,1}(M)$
the Maurer-Cartan equation
\begin{equation}\label{_MC_for_1,1_Equation_}
\bar \6\gamma_p = \sum_{i+j=p-1} [\gamma_i, \gamma_j], p=1,2, 3, ...
\end{equation}
has a solution
$(\gamma_1, \gamma_2, ..., )$, with $\gamma_i \in \Lambda^{1,1}(M)$
$\6$-closed. Here the commutator $[\gamma_i, \gamma_j]$
is understood as the Schouten bracket, see \ref{_hol_symp_TT_Theorem_}.

\hfill

Unobstructedness of deformations of 
$(M,I, \Omega)$ implies that any $\bar\6$- and $\6$-closed form $\gamma_0\in \Lambda^{1,1}(M)$
is tangent to a one-parametric family of deformations 
of holomorphic symplectic structures. 
Using the inverse function theorem, we obtain that the
deformation space of holomorphic symplectic structures, if unobstructed, 
is locally biholomorphic to an open ball in $H^1(M, \Ham_M)$
(see \cite{_KV:deformations_} for a formal proof). 

\hfill

%

The main result of this section is the following theorem,
generalizing the Bogomolov theorem on unobstructedness
of deformations of hyperk\"ahler manifolds (\cite{_Bogomolov:TT_}).

\hfill

\theorem \label{_TT_for_HS_Theorem_}
Let $(M,I, \Omega)$ be a compact holomorphically symplectic
manifold (not necessarily K\"ahler). Assume that the Dolbeault
cohomology group $H^{0,2}_{\bar \6}(M)= H^2(\calo_M)$ is generated
by $\6$-closed $(0,2)$-forms, and all $\6$-exact holomorphic 3-forms on $M$ vanish.  Then 
\begin{description}
\item[(a)]
the holomorphic symplectic
deformations of $(M,I, \Omega)$ are unobstructed.
\item[(b)]
If, in addition, all classes in the Dolbeault
cohomology group $H^{1,1}_{\bar \6}(M) $ are represented by closed $(1,1)$-forms,
then the complex deformations of $M$ are also unobstructed,
and all sufficiently small complex deformations remain
holomorphically symplectic.
\end{description}

\hfill

\proof (See also \cite[Theorem 2.2]{_KV:deformations_}).
Let $\alpha, \beta \in \Lambda^{1,*}(M)$
be $\6$-closed forms, and $[\alpha, \beta]$
the Schouten bracket on $\Lambda^{1,*}(M)$
defined as in \ref{_hol_symp_TT_Theorem_}.
Then \ref{_hol_symp_TT_Theorem_}
gives 
\begin{equation}
[\alpha, \beta]= \6\Lambda_\Omega(\alpha\wedge \beta).
\end{equation}
Suppose we have solved the Maurer-Cartan equation
\eqref{_MC_for_1,1_Equation_} for all $p < n$,
and all $\gamma_p$ with $p<n$ are $\6$-closed.
Then \eqref{_MC_for_1,1_Equation_} becomes
\begin{equation}\label{_MC_via_Lambda_Equation_}
\bar\6\gamma_n= \sum_{i+j=n-1} \6\Lambda_\Omega(\gamma_i\wedge \gamma_j).
\end{equation}
By \eqref{_triple_commu_Equation_}, the right hand side of
\eqref{_MC_via_Lambda_Equation_} is $\bar\6$-closed.
It is $\bar\6$-exact by the following simple lemma,
applied to $\rho=\Lambda_\Omega(\gamma_i\wedge \gamma_j)$.

\hfill

\lemma
In assumptions of \ref{_TT_for_HS_Theorem_},
let $\rho\in \Lambda^{0,2}(M)$ be a form which satisfies
$\bar\6\6\rho=0$. Then $\6\rho$ is $\6\bar\6$-exact.

\hfill

\proof 
Since $\bar \6\6\bar \rho=0$,
the $(3,0)$-form $\6\bar \rho$ is $\6$-exact and holomorphic.
By assumptions of \ref{_TT_for_HS_Theorem_}, it vanishes.
 Then
$\bar\6\rho=0$. Since all $\bar\6$-cohomology
classes in $\Lambda^{0,2}(M)$ can be represented
by closed forms, there exists
a closed $(0,2)$-form $\rho'=\rho-\bar\6\mu$.
This gives $\6\rho=\6\bar\6\mu$.
\endproof

\hfill

We have proven \ref{_TT_for_HS_Theorem_} (a).
To prove (b), we start with $\eta\in H^{1,1}_{\bar\6}(M)= H^1(TM)$,
take a $\6, \bar\6$-closed representative, and solve the
Maurer-Cartan equation \eqref{_MC_for_1,1_Equation_}
as above. This gives a deformation of complex structures,
which is by construction holomorphically symplectic
(\ref{_closed_gamma_holo_symp_Claim_}).
\endproof

\hfill

Let $U$ be the Kuranishi deformation space of 
the holomorphically symplectic manifold $M$, that is
of pair $(I,\Omega)$ of complex and holomorphically
symplectic structure.
Consider the period map $\Per:\; U \arrow H^2(M, \C)$
mapping $(I,\Omega) \in U$ to the the cohomology class 
of the holomorphically symplectic form $\Omega$.
In \ref{_closed_gamma_holo_symp_Claim_}, we interpret
the holomorphically symplectic deformations as 
closed $(1,1)$-forms. It is not hard to see that 
the differential $d\Per$ of this map takes $\gamma_0 \in H^{1,1}(M,I)= T_I U$
and maps it to the same vector in $T_{[\Omega]} H^2(M, \C)= H^2(M, \C)$
(see also \cite{_KV:deformations_}).

\hfill

\remark\label{_Omega_defo_Remark_}
Let $S\subset H^2(M, \C)$ be the set of all cohomology classes of form $\Omega$,
where $\Omega$ is the cohomology class of a holomorphically symplectic form
for some $I\in U$. From the above argument it is clear $T_{[\Omega]} S$ contains
$(1+a) \Omega +  S_0$, where $S_0$ is a neighbourhood
of 0 in $H^{1,1}(M, \C)$, and $a\in \C$ sufficiently small
number. This implies the following useful corollary.

\hfill

\corollary \label{_Grassmann_periods_Corollary_}
In assumptions of \ref{_TT_for_HS_Theorem_}, 
let $U$ be the Kuranishi deformation space
of holomorphically symplectic structures satisfying
assumptions of \ref{_TT_for_HS_Theorem_}. Let $\Gr(2, H^2(M, \R))$
be the  Grassmann space of 2-planes in $H^2(M, \R)$.
Consider the period map
$\Per:\; U \arrow \Gr(2, H^2(M, \R))$
taking $(I, \Omega) \in U$ and mapping it to the plane generated by
$\langle \Re\Omega, \Im \Omega\rangle\subset H^2(M, \R)$. Then 
$\Per(U)$ is open in $\Gr(2, H^2(M, \R))$.

\hfill

\proof By construction, the image of $(I,\Omega)$ under the period map is
$\langle \Re\Omega, \Im \Omega\rangle\subset H^2(M,
\R)$. From \ref{_Omega_defo_Remark_} 
it follows that $T_{[\Omega]} \Per(U)$
contains $\langle (1+a) \Re\Omega + S_{\Re}, (1+b)\Im
\Omega + S_{\Im} \rangle\subset H^2(M,\R)$, where $S_\Re, S_\Im$ are  neighbourhoods
of 0 in $H^{1,1}(M, \R)$, corresponding to 
the real and imaginary part of $\alpha\in S_0$,
and $a, b \in \R$ sufficiently small real numbers.   
By \ref{_TT_for_HS_Theorem_} it
follows that the deformation space is smooth, impliying
that $\Per(U)$ is open in $\Gr(2, H^2(M, \R))$.
\endproof


%


\section{Beauville-Bogomolov-Fujiki form} \label{BBF-form-section}


In this section we generalize the construction
of Bogomolov-Beauville-Fujiki (BBF) form to non-K\"ahler
holomorphically symplectic manifolds. 

\subsection{BBF form on non-K\"ahler holomorphically symplectic manifolds}

\definition\label{_Hodge_deco_Definition_}
Let $M$ be a complex manifold.
We say that {\bf $H^2(M)$ admits the Hodge decomposition}
if any cohomology class in $H^2(M)$ can be represented by a sum of
closed $(p,q)$-forms.

\hfill

\example
All K\"ahler manifolds admit the Hodge decomposition (this is how
the Hodge decomposition was originally discovered).
Moreover, all compact complex surfaces also admit the Hodge
decomposition in $H^2(M)$ (\cite[Theorem IV.2.8]{_Barth_Peters_Van_de_Ven_}).

\hfill

This is the main
result of this paper.

\hfill

\theorem\label{_BB_main_Theorem_}
Let $M$ be a compact holomorphically symplectic manifold,
$\dim_\C M=2n$, admitting the Hodge decomposition on $H^2(M)$.
Assume that all $\6$-exact holomorphic 3-forms on $M$ vanish,
and $H^{0,2}(M)=H^{2,0}(M)=\C$. Then the space
$H^2(M)$ is equipped with a bilinear symmetric form
$q$ such that for any $\eta\in H^2(M)$, one has
$\int_M \eta^{2n} = \lambda q(\eta, \eta)^n$, where
$\lambda$ is a fixed constant.

\hfill

\proof
Let $U$ be the local deformation space of 
holomorphic symplectic structures on $M$.
Shrinking $U$ if necessary, we may assume that
all complex structures $I\in U$ satisfy
assumptions of \ref{_BB_main_Theorem_}.

Consider the Hodge decomposition on
$H^2(M)$ induced by the complex structure $I\in U$.
This gives a $U(1)$-action $\rho_I$
on $H^2(M)$, with $\rho_I(t)$ acting as $e^{2\pi\1(p-q) t}$ on $H^{p,q}(M)$.
Clearly, the polynomial $Q(\eta):=\int_M \eta^{2n}$ is $\rho_I$-invariant.

Now, \ref{_BB_main_Theorem_}
is implied by the following algebraic proposition, applied
to the group $G$ generated by all $\rho_I$, $I\in U$, and
$V=H^2(M, \R)$.

\hfill

\proposition\label{_BBF_algebraic_via_Hodge_deco_Proposition_}
Let $V$ be a real vector space equipped with an action of a Lie group $G$,
and $Q$ a $G$-invariant polynomial function. Let $S\subset \Gr(2, V)$
be an open subset in the Grassmannian of 2-planes. Assume that 
for any $P\in S$, there exists a subgroup $\rho_P\subset G$ isomorphic
to $S^1$ acting by rotations on $P$ and trivially on $V/P$.
Then $Q$ is proportional to $q^n$, where $q$ is a quadratic form 
on $V$.

\hfill

We prove \ref{_BBF_algebraic_via_Hodge_deco_Proposition_} below. To deduce
\ref{_BB_main_Theorem_} from 
\ref{_BBF_algebraic_via_Hodge_deco_Proposition_},
we notice that each $I\in U$ is associated with the
Hodge rotation $\rho_I$ as above, and the set of 
2-planes generated by $\langle \Re\Omega, \Im \Omega\rangle$
for all $\Omega\in H^{2,0}(M,I)$ is open in $\Gr(2, H^2(M, \R))$
(\ref{_Grassmann_periods_Corollary_}).
Then we apply \ref{_BBF_algebraic_via_Hodge_deco_Proposition_}
to obtain that the polynomial $Q(\eta):=\int_M \eta^{2n}$
is proportional to $q^n$, where $q$ is a quadratic form.
This proves \ref{_BB_main_Theorem_}.
\endproof

\subsection{Polynomial invariants of Lie groups}
\label{_Poly_inva_Subsection_}

In this subsection, we prove \ref{_BBF_algebraic_via_Hodge_deco_Proposition_}.

\hfill

\noindent {\bf \ref{_BBF_algebraic_via_Hodge_deco_Proposition_}:\ }
Let $V$ be a real vector space equipped with an action of a Lie group $G$,
and $Q$ a $G$-invariant polynomial function. Let $S\subset \Gr(2, V)$
be an open subset in the Grassmannian of 2-planes. Assume that 
for any $P\in S$, there exists a subgroup $\rho_P\subset G$ isomorphic
to $S^1$ acting by rotations on $P$ and trivially on $V/P$.
Then $Q$ is proportional to $q^n$, where $q$ is a quadratic form 
on $V$.

\hfill

\proof
Let $P\in S$ be a 2-plane in $V$. Any rotation-invariant polynomial
function on $\R^2$ is a power of quadratic form, hence 
$Q\restrict P= \lambda q^n \restrict P$. When $n$ is odd, 
the $n$-th root of $Q$ is well defined. When $n$ is even,
the restriction $Q\restrict P$ does not change sign,
hence $Q$ does not change sign on the set $U_S\subset V$
of all vectors passing through planes $P\in S$.
The function $q:=\sqrt[n]{\pm Q}$
is well defined  on the whole of $V$ when $n$ is odd,
and on $U_S$ when it is even. Also, this function is polynomial of second degree on
all hyperplanes $P\in S$. Therefore, the second derivative
$\frac{d^2}{dxdy}q$ vanishes on $V$ or $U_S$
when $\langle x, y\rangle \in S$. This function
is real algebraic, hence by the analytic
continuation principle, $\frac{d^2}{dxdy}q=0$
everywhere, and this implies that $q$ is a 
second degree polynomial.

This proves  \ref{_BBF_algebraic_via_Hodge_deco_Proposition_}.
\endproof


\section{Bogomolov-Guan manifold and its deformations}
\label{_BG_Section_}


\subsection{Cohomology of Hilbert scheme of complex surfaces} \label{coh-Hilb-scheme}

\definition
Let $M$ be a compact complex manifold.
We say that the {\bf Hodge decomposition holds for $H^p(M)$}
if any cohomology class $x\in H^p(M)$ can be represented by a form
$\xi\in \Lambda^{p}(M)$ which satisfies
$\6\xi=\bar\6\xi=0$. For an equivalent definition, 
see \ref{_Hodge_deco_Definition_}.

\hfill

\proposition \label{Hodge-decomp-Hilbert}
Let $M$ be a compact complex surface, and $\Hilb^n(M)$
the $n$-th Hilbert scheme of points on $M$. Then $\Hilb^n(M)$
admits the Hodge decomposition for $H^2(M)$.

\hfill

\proof As follows from \cite[Theorem IV.2.8]{_Barth_Peters_Van_de_Ven_},
$M$ admits the Hodge decomposition for $H^2(M)$. Then the Hodge decomposition in $H^2$ for $M^{[n]}$ follows from the fact that $H^2(M^{[n]})=H^2(M) \oplus [E]$, where $E$ is the exceptional divisor, and $[E]$ is fundamental class.
 \endproof

\hfill

\begin{remark}
We can also use \cite[Lemma 2.10]{_LQZ_Hilb}, where 
the generators of $H^2(N^{[n]})$ are described
using the results of G\"ottsche and Nakajima.
\end{remark}

\hfill

\claim \label{holomorphic-Sym}
Let $M$ be a compact complex surface, $M^{[n]}$
the $n$-th Hilbert scheme of points on $M$, and $M^{(n)}$ its $n$-th symmetric
power. Then $H^0(\Omega^i(M^{[n]}))= H^0(\Omega^i(M^{(n)}))$.

\hfill

\proof 
Birational maps induce an isomorphism on global holomorphic $k$-forms for any $k$.
\endproof

\subsection{Bogomolov-Guan manifold and its deformations}

In this section we consider the Bogomolov-Guan manifold $Q$
constructed in Subsection \ref{BG-construction_section}. Recall that
$F^{[n]}$ is the preimage of $0 \in E$ of the map 
$\pi:\; S^{[n]} \rightarrow E$. Let $W$ be the leaf space of
characterstic foliation on $F^{[n]}$. Then Bogomolov-Guan
manifold $Q$ is a degree $n^2$ covering of $W$. By
construction we have the following diagram (part of
\eqref{BG_construction})
\begin{displaymath} \label{BG_constr-cohomology}
\begin{diagram}[labelstyle=\tiny]
 &  & H^{[n]} & \rTo & F^{[n]}\\
 & &\dTo^{S^1}& & \dTo_{E_L} \\
 Q &\rTo^{p_r} & R&\rTo^{p}  & W, \\ \end{diagram}
\quad
\end{displaymath}
where $H^{[n]}$ is the preimage of zero of torus $T^3$, and the map $p_w:\; Q \rightarrow W$ is $p_w=p \circ p_r$.




\hfill

\begin{definition}
An {\bf orbispace} is a topological space $M$, 
equipped with:
\begin{enumerate}
	\item  A structure of a groupoid (the points of $M$
are objects of a groupoid category);
\item  A covering $\{U_i\}$;
 \item Continuous maps $\phi_i:\; V_i\arrow U_i$, where each
$V_i$ is equipped with a rigid action of a finite group $G_i$,
satisfying the following properties: 
\begin{enumerate}
\item $\phi_i:\; V_i \arrow V_i/G_i=U_i$ is the quotient map.
\item For each $x\in M$ and $U_i\ni x$, 
the group $\Mor(x,x)$ (i.e. the monodromy group of an orbipoint) is equal to the stabilizer $\St_G(x)$ of $x$ in $G_i$.
\end{enumerate}
\end{enumerate}
\end{definition}

\hfill

\begin{definition}
 A {\bf smooth orbifold} is an orbifold $M$ equipped with
a sheaf of functions $ C^\infty (M)$ in such a way that for each 
$U_i=V_i/G_i$, the corresponding ring of sections
$ C^\infty (U_i)$ is identified with the  
ring of $G_i$-invariant smooth functions on $V_i$.
\end{definition}

\hfill

Cohomology of a smooth orbifold are cohomology of its de Rham algebra.

\hfill






\proposition \label{holomorphic-hilbert}
Let $S^{[n]}$ be a Hilbert scheme of points of
Kodaira-Thurston manifold. Then all holomorphic 3-forms on
$S^{[n]}$ vanish. Moreover, the ring of holomorphic differential forms on
$S^{[n]}$ is isomorphic to
$\frac{\C[\theta, \Omega]}{\langle \theta^2=0,
  \Omega^{n+1}=0 \rangle}$,
where $\Omega$ is a holomorphic symplectic form, and
$\theta$ is a 1-form.

\hfill

\proof
Notice that the blow-up $S^{[n]}\arrow S^{(n)}$
is bimeromorphic and induces an isomorphism on
the global holomorphic differential forms. 
The ring of holomorphic forms on a Kodaira-Thurston surface
is 
\[ \frac{\C[\theta, \Omega]} {\langle \theta^2=\Omega^2=\theta
\wedge \Omega =0 \rangle},
\] where $\theta$ has degree one, and $\Omega$ has degree
two. Holomorphic forms on the Hilbert scheme $S^{[n]}$
are identified with $H^{*,0}(S^n)^{\Sigma_n}$,
where $\Sigma_n$ is the symmetric group acting on $S^n$
and $H^{*,0}(S^n)^{\Sigma_n}$ denotes the
$\Sigma_n$-invariant
holomorphic forms.

Denote by $\pi_i:\; S^n \arrow S$
the projection of $S^n$ to its $i$-th component.  Let
$\theta_i=\pi^*_i\theta$, and $\Omega_i=\pi_i^*\Omega$. 
Then $H^{*,0}(S^n)^{\Sigma_n}=\C[\Omega_i, \theta_j]^{\Sigma_n}$ is the space of
$\Sigma_n$-invariant polynomials of $\theta_i$ and $\Omega_j$.
Any monomial $\Omega_{i_1} \Omega_{i_2}... \Omega_{i_k}
\theta_{i_{k+1}} ... \theta_{i_{k+m}}$ vanishes if
some of the indices coincide, because $\Omega\wedge\theta$,
$\theta\wedge \theta$ and $\Omega\wedge \Omega$
vanish on $S$. Symmetrizing such a monomial with
$\Sigma_n$, we obtain a sum
\[
R_{k,m}:=\sum_{\sigma\in \Sigma_n}\Omega_{i_{\sigma_1}} \Omega_{i_{\sigma_2}}... \Omega_{i_{\sigma_k}}
\theta_{i_{\sigma_{k+1}}} ... \theta_{i_{\sigma_{k+m}}},
\]
where all indices are pairwise distinct.
As shown above, the space of holomorphic differential forms
on $S^{(n)}$ and on $S^{[n]}$ is generated as a vector
space by $R_{k,m}$.

Consider now a holomorphic form
\[ R(k,m):=\left(\sum_{i=1}^n \Omega_i\right)^k \left(\sum_{i=1}^n
\theta_i\right)^{m}.
\]
Decomposing $R(k,m)$ onto a sum of monomials as above,
we find that the only non-zero summands are monomials
with nonequal indices. After averaging with $\Sigma_n$,
each of these monomials gives $R_{k,m}$.
This gives $ R(k,m)=\const R_{k,m}$. We have shown that
the ring of holomorphic differential forms on 
$S^{(n)}$ and on $S^{[n]}$ is generated as a ring by 
$\tilde \Omega:=\sum_{i=1}^n \Omega_i$ and $\tilde
\theta:=\sum_{i=1}^n\theta_i$. 
The relations $\tilde \theta^2=0$
and $\tilde \Omega^{n+1}=0$ are apparent because 
$\dim_\C S^{[n]}=2n$ and $\tilde \theta$ is a 1-form. However,
$\tilde \Omega^{m}\neq 0$ and $\tilde \Omega^{m}\wedge \tilde \theta \neq 0$  
for $m\leq n$ because $\tilde \Omega$ is non-degenerate.
Therefore, the ring of holomorphic differential forms on
$S^{[n]}$ is isomorphic to
$\frac{\C[\tilde \theta, \tilde \Omega]}{\langle \tilde \theta^2=0,
  \tilde \Omega^{n+1}=0 \rangle}$. 
\endproof

\hfill

\proposition \label{cohomology-F} 
Consider the complex
manifold $F^{[n]}$ constructed above as the preimage of a
map $\pi:\; S^{[n]} \rightarrow E$. Then $H^0(F^{[n]},
\Omega^i_{F^{[n]}}) = \left\langle \Omega^i
\right\rangle$, where $\Omega$ is the holomorphic
symplectic form on $S^{[n]}$.

\hfill
 
\proof
We start by proving that 
\[ H^0(S^{[n]}, \Omega^*_{S^{[n]}})=H^0 \left( F^{[n]},
\Omega^*_{S^{[n]}}\restrict {F^{[n]}} \right).
\]
Indeed, $F^{[n]}$ is the fiber of a locally trivial
complex fibration $S^{[n]} \rightarrow E$. Consider the
following exact sequence
\begin{multline*}
 0 = H^0\left(S^{[n]}, \Omega^i S^{[n]}\left(-F^{[n]}\right)\right) \rightarrow H^0\left(S^{[n]}, \Omega^i S^{[n]}\right) \rightarrow \\
 \rightarrow H^0\left(F^{[n]}, \Omega^i S^{[n]}\restrict {F^{[n]}}\right) \rightarrow H^1\left(S^{[n]}, \Omega^i S^{[n]}\left(-F^{[n]}\right)\right) =0,
\end{multline*} 
where for the left and right equations we use
$\rho_*S^{[n]}\left(-F^{[n]}\right) = \bigoplus \mathcal{O}(-1) = 0$
on $E$, since $E$ is the elliptic curve. 

Consider the isotrivial fibration
$\pi:\;S^{[n]}\arrow E$, with fiber $F^{[n]}$. Since this
fibration is isotrivial, one has a sheaf decomposition
$\Omega^1_{S^{[n]}}\restrict {F^{[n]}}=
  \Omega^1_{F^{[n]}}\oplus N^*{F^{[n]}}$,
where $N^*{F^{[n]}}=\calo_{F^{[n]}}$ is the (trivial) conormal
    bundle to the fiber of $\pi$. Then
\[
H^0\left(F^{[n]}, \Omega^i S^{[n]}\restrict
{F^{[n]}}\right)= H^0\left(\Omega^*_{F^{[n]}}\right)
\oplus  \theta \wedge H^0\left(\Omega^*_{F^{[n]}}\right),
\]
where $\theta$ denotes the generator of the
conormal component of $\Omega^1_{S^{[n]}}\restrict{F^{[n]}}$. 
In the proof of \ref{holomorphic-hilbert} this form was
denoted as $\tilde \theta= \sum_i \theta_i$.
Now, \ref{cohomology-F} would follow if
we prove that the restriction $\theta
\restrict{F^{[n]}}$ vanishes. This is clear, because
$\theta$ is a pullback of a holomorphic form on $E$,
hence restricts trivially on the fibers of $\pi:\;S^{[n]}\arrow E$.
\endproof 

\hfill

\proposition \label{_BG-conditions-BB-existence_Proposition_}
Let $Q$ be a Bogomolov-Guan manifold (\ref{_Bogomolov_Guan_Definition_}).
Then $Q$ admits the Hodge decomposition in $H^2(M)$
and all holomorphic 3-forms on $Q$ vanish.

\hfill

\proof 
Existence of the Hodge decomposition follows from \cite[Theorem 2]{Gu3}; see also Subsection \ref{coh-Hilb-scheme}.

Now we will prove the vanishing of holomorphic 3-forms on $Q$.

\hfill

{\bf Step 1.} By \ref{cohomology-F} this is true for $F^{[n]}$. 

\hfill

{\bf Step 2.} We prove that 
\[ H^0(W, \Omega^i_W) =
\left\langle \Omega^i \right\rangle,
\] 
where $W$ is the leaf space of characteristic foliation,
considered as a holomorphically symplectic orbifold,
and $\Omega$ its holomorphic symplectic form.
Using the dominant holomorphic map $ F^{[n]}\arrow W$,
we can realize any holomorphic form on $W$ as a 
holomorphic form on $F^{[n]}$. This gives an embedding
$H^0(\Omega^i W) \subset H^0(\Omega^i F^{[n]})$. Then  \ref{cohomology-F}
implies that $H^0(W, \Omega^i_W)$ is generated by the
powers of its holomorphic symplectic form.

\hfill

{\bf Step 3.} Consider the unramified orbifold map
$p_w\colon Q \rightarrow W$, which is obtained as a
composition of two maps of degree $n$, and $\Gamma$ is a cyclic
deck transform group. Then we have
$H^0(\Omega^*Q)=H^0(\Omega^*W \otimes
p_{w*}\mathcal{O}_Q)$, which is seen as follows.

Differential forms on Bogomolov-Guan manifold are obtained as
$\Omega^* Q = p_w^* \Omega^* W$, because
the map $p_w:\;Q \arrow W$ is the universal cover in the orbifold category.
Applying the pushforward,
we obtain $p_{w*}\Omega^*Q=p_{w*}p_w^* \Omega^*
W=\Omega^*W \otimes p_{w*}\mathcal{O}_Q$. The latter term
decomposes as $p_{w*}\mathcal{O}_Q = \bigoplus_{\chi \in
  \Gamma} L_\chi$, where $\chi$ are characters of the
cyclic group $\Gamma$, and $L_i$ are the sheaves of $\chi$-automorphic
functions.

\hfill

{\bf Step 4.} To prove that $H^0(\Omega^*Q)=H^0(\Omega^*W \otimes
p_{w*}\mathcal{O}_Q)$ it remains to show that
$H^0(\Omega^*W\otimes L_\chi)=0$ if $\chi$ is
non-trivial. For a trivial character $\chi_0$, the line bundle $L_{\chi_0}$
is trivial. Then $H^0(\Omega^*W)=H^0(\Omega^*W\otimes L_{\chi_0})$, 
and  we can apply Step 2.

Consider the homotopy exact sequence associated with the
isotrivial fibration $\pi\colon S^{[n]} \rightarrow E$:
\[0=\pi_2(E) \rightarrow \pi_1(F^{[n]}) \rightarrow \pi_1(S^{[n]}) \rightarrow \pi_1(E) \rightarrow 0
\] 
Then we have $\tors(\pi_1(F^{[n]}) \hookrightarrow \tors(\pi_1(S^{[n]})$, where $\tors$ is the torsion
subgroup. 
The group $\pi_1(E)=\Z^2$ is a projective
$\Z$-module. Therefore, $\tors(\pi_1(F^{[n]}))$ is a direct
summand in $(\pi_1(S^{[n]})$. Hence, $\lambda^*L_\chi$ can
be extended on $S^{[n]}$. Abusing this notation,
we use the same letter $L_\chi$ for this line bundle and its
extension to $S^{[n]}$. Then
\[
H^0(\Omega^*F^{[n]} \otimes L_\chi) \hookrightarrow
H^0(\Omega^* S^{[n]} \otimes L_\chi)=0,
\]
One has $H^0(\Omega^* S^{[n]} \otimes L_\chi)=0$ for any
$\chi\neq 0$.
Indeed, the appropriate $|\Gamma|$-cover of $S^{[n]}$ is
the Hilbert scheme $S'^{[n]}$ of a Kodaira surface $S'$,
which is a $|\Gamma|$-cover of $S$, and
$H^0(\Omega^*S'^{[n]})=\left\langle \Omega^i
\right\rangle$. Using the triviality of $\Gamma$-action on
$H^0(\Omega^*S'^{[n]})$, we obtain vanishing of
$H^0(\Omega^*W\otimes L_\chi)$ if $\chi$ is non-trivial. 

\hfill

{\bf Step 5.} Using the dominant map
$F^{[n]}\arrow W$
as in Step 2, we obtain an embedding
 $H^0(\Omega^*_W\otimes L_\chi) \subset H^0(\Omega^*_{F^{[n]}}\otimes L_\chi)$ 
for any character $\chi$. When $\chi$ is non-trivial,
one has $H^0(\Omega^*_{F^{[n]}}\otimes L_\chi)=0$ 
because all holomorphic forms on
$S'^{[n]}$ are $\Gamma$-invariant.

\hfill

Combining all steps we have 
\[ H^0(\Omega^*Q) \xlongequal{\text{Step 3}}
\bigoplus_{\chi \in \Gamma} H^0(\Omega^*W\otimes L_\chi)\xlongequal{\text{Step
    4-5}}H^0(W, \Omega^i)\xlongequal{\text{Step 2}}
\left\langle \Omega^i \right\rangle.
\]
\ref{_BG-conditions-BB-existence_Proposition_} is proven. \endproof

\hfill

\corollary \label{BBF-form-BG}
The deformation space of a Bogomolov-Guan manifold $Q$
is smooth in a neighbourhood of $Q$, and all sufficiently
small deformations are holomorphically symplectic. Moreover,
$H^2(Q)$ is equipped with a bilinear symmetric form
$q$ (Beauville-Bogomolov-Fujiki form) such that for any $\eta\in H^2(Q)$, one has
$\int_Q \eta^{2n} =  q(\eta, \eta)^n$,
where $2n=\dim_\C Q$.

\hfill

\proof The conditions required by \ref{_TT_for_HS_Theorem_}
follow from \ref{_BG-conditions-BB-existence_Proposition_}. This implies the
smoothness of the deformation space. The existence of the
Beauville-Bogomolov-Fujiki form on a Bogomolov-Guan manifold $Q$ follows from \ref{_BB_main_Theorem_}. \endproof

\hfill

\conjecture Let $Q$ be a Bogomolov-Guan manifold and $q$
the Beauville-Bogomolov-Fujiki form on $H^2(Q,\R)$. Then $q$ is non-degenerate.

\hfill

{\bf Acknowledgements:} We would like to thank
F. Bogomolov,  D. Kaledin and G. Papayanov for useful comments.

\footnotesize
{

\noindent
{\sc Nikon Kurnosov\\
Department of Mathematics,\\
University of Georgia,\\
Athens, GA, USA, 30602\\
also:\\
\sc  Laboratory of Algebraic Geometry,\\
National Research University HSE,\\
Department of Mathematics,
6 Usacheva Str. Moscow, Russia\\
\tt nikon.kurnosov@gmail.com}\\

\hfill

\noindent
{\sc Misha Verbitsky\\
 Instituto Nacional de Matem\'atica Pura e
              Aplicada (IMPA) \\ Estrada Dona Castorina, 110\\
Jardim Bot\^anico, CEP 22460-320\\
Rio de Janeiro, RJ - Brasil \\
also:\\
\sc Laboratory of Algebraic Geometry,\\
National Research University HSE,\\
Department of Mathematics, 6 Usacheva Str. Moscow, Russia\\
\tt  verbit@impa.br}.
}


{\small 
\begin{thebibliography}{HKLR}

\bibitem[AFM]{_AFM:KT_}
L. C. de Andr\'es, M. Fern\'andez, J. Menc\`ia, {\em Curvature and complex geometry on the
  Kodaira-Thurston manifold}, Proceedings of the Workshop
on Curvature Geometry (Lancaster, 1989), 95-105, ULDM
Publ., Lancaster, 1989. 

\bibitem[A]{_Angela_}
D. Angella, {\em Cohomological Aspects in Complex Non-K\"ahler Geometry}, Lecture Notes in Mathematics 2095, Springer, 2014.

\bibitem[BT]{_Babenko_Taimanov_} 
I. K. Babenko, I. A. Taimanov,
{\em Massey products in symplectic manifolds},
math.SG/9911132, \textit{Sb. Math.}, {\bf 191} (2000), pp. 1107--1146.

\bibitem[BL]{_Bakker_Lehn_}
B. Bakker, C. Lehn, {\em The global moduli theory of symplectic varieties}, ArXiv, math.AG: 1812.09748.

\bibitem[BK]{_Barannikov_Kontsevich_} 
S. Barannikov, M. Kontsevich, 
{\it Frobenius Manifolds and Formality 
of Lie Algebras of Polyvector Fields},
 \textit{Int. Math. Res. Not.}, {\bf 1998}, no. 4, pp. 201--215.



\bibitem[BHPV]{_Barth_Peters_Van_de_Ven_}
W. Barth, K. Hulek, C. Peters, A. Van de Ven, 
{ Compact complex surfaces},  Springer Verlag, 2004.


\bibitem[Bea]{_Beauville_} 
A. Beauville, {\em 
Varietes K\"ahleriennes dont la premi\`ere classe de Chern est
nulle.},  \textit{J. Diff. Geom.}, {\bf 18} (1983), pp. 755 -- 782.






\bibitem[Bo1]{_Bogomolov:decompo_}  
F. A. Bogomolov, {\em On the decomposition of 
K\"ahler manifolds with trivial canonical class}, \textit{Math. USSR-Sb.},
{\bf 22} (1974), pp. 580 -- 583.




\bibitem[Bo3]{B1} F. Bogomolov, {\em On Guan's examples of simply connected non-K\"ahler compact
complex manifolds}, \textit{Amer. Journ. of Math.}, \textbf{118}, Number 5 (1996), pp. 1037--1046.

\bibitem[Bo4]{_Bogomolov:TT_}
F. Bogomolov, {\em Hamiltonian K\"ahler manifolds}, Dokl. Akad. Nauk SSSR 243 (1978), 1101--1104; Soviet Math. Dokl., 19, (1979), 1462--1465.








\bibitem[G]{G} P. Gauduchon, {\em La 1-forme de torsion d'une variete hermitienne compacte}, \textit{Math.Ann.}, 1984.

\bibitem[Gu1]{Gu1} D. Guan, {\em Examples of compact holomorphic symplectic manifolds which admit no K\"ahler structure},
Geometry and Analisys on Complex Manifolds—Festschrift for Professor Kobayashi S. 60th
Birthday, World Scientific, Teaneck, NJ, 1994, pp. 63–74.

\bibitem[Gu2]{Gu2} D. Guan, { \em Examples of compact holomorphic symplectic manifolds which are not Kahlerian II}, \textit{Invent. math.}, \textbf{121.1} (1995), pp. 135--146. 

\bibitem[Gu3]{Gu3} D. Guan, {\em Examples of compact holomorphic symplectic manifolds which are not Kahlerian III}, \textit{Int. J. Math.}, \textbf{06}, 5, pp. 709 -- 718 (1995). 


\bibitem[C]{_Calabi_} 
E. Calabi,
{\em Metriques k\"ahleriennes et fibr\`es holomorphes}, 
\textit{Ann. Ecol. Norm. Sup.}, {\bf 12} (1979), pp. 269--294.  


\bibitem[F]{_Fujiki:HK_}  
A. Fujiki {\em On the de Rham Cohomology Group of a Compact 
K\"ahler Symplectic Manifold}, \textit{Adv. Stud.
Pure Math.}, \textbf{10} (1987), pp. 105--165.



\bibitem[G]{_Ghys_}
\'E. Ghys, {\em D\'eformations des structures complexes sur les
 espaces homog\`enes de SL(2,C),}
\textit{J. Reine Angew. Math.}, {\bf 468} (1995), 113--138.

\bibitem[I]{_Iacono_}
D. Iacono, {\em On the abstract Bogomolov-Tian-Todorov Theorem}, \textit{Rend. Mat. Appl.} (7).
Volume \textbf{38}, (2017), pp. 175 -- 198.

\bibitem[KV]{_KV:deformations_}
D. Kaledin, M. Verbitsky
{\em Period map for non-compact holomorphically symplectic manifolds},
GAFA, \textbf{12} (2002), no. 6, pp. 1265--1295.

\bibitem[Ki]{_Kirshner}
T. Kirschner, \textit{Period  mappings  with  applications  to  symplectic  complex  spaces},  v. 2140 of \textit{Lecture Notes in Mathematics}, Springer, Cham, 2015.

\bibitem[Kod]{_Kodaira:surfaces_1_}
Kodaira, K.
{\em On the structure of compact complex analytic surfaces. I,}
Amer. J. Math. 86 (1964), 751-798. 

\bibitem[Kon]{_Kontsevich:lectures_}
M. Kontsevich, {\em Topics in algebra: deformation theory},
notes by Alan Weinstein,
\url{http://www1.mat.uniroma1.it/people/manetti/DT2011/Kontsevich.pdf}

\bibitem[L]{_de_Leon:onThurston_}
M. de Le\'on, {\em 
Sur une conjecture de Thurston,} C. R. Acad. Sci. Paris S'er. I Math. 301 (1985), no. 16, 771. 

\bibitem[LQZ]{_LQZ_Hilb} W. Li, Z. Qin, Q. Zhang, {\em On the geometry of the Hilbert schemes of points in the projective plane}, ArXiv: 0105213.



\bibitem[Na1]{_Namikawa_} Y. Namikawa, {\em On deformations of Q-factorial symplectic varieties}, \textit{J.
Reine Angew. Math. (Crelle Journ.)}, \textbf{599} (2006), pp. 97--110.

\bibitem[Na2]{_Namikawa_form}
Y. Namikawa. {\em Extension of 2-forms and symplectic varieties}, \textit{J. Reine Angew. Math.}, \textbf{539} (2001), pp. 123--147.


\bibitem[Th]{_Thurston:Kodaira_}
W. Thurston, {\em Some simple examples of symplectic manifolds}, Proc. Amer. Math.
Soc. 55 (1976) 467-468.


\bibitem[Ti]{_Tian:TT_} 
G. Tian, {\em Smoothness of the universal
deformation space of compact Calabi-Yau manifolds and its
Petersson-Weil metric}, in {\em Math. Aspects of String Theory},
S.-T. Yau, ed., Worlds Scientific, 1987, pp. 629--646.

\bibitem[To1]{_Todorov:MPIM_}
A. Todorov
{\em	Every holomorphic symplectic manifold admits a K\"ahler metric},
MPIM preprint 1985-43, \url{https://www.mpim-bonn.mpg.de/preblob/5418}.

\bibitem[To2]{_Todorov:TT_} 
A. Todorov, {\em The Weil-Petersson geometry of
the moduli space of $SU(n \geq 3)$ (Calabi-Yau) manifolds}, \textit{Comm.
Math. Phys.}, {\bf 126} (1989), pp. 325--346.




%
%



\bibitem[V1]{_V:Mirror_}
M. Verbitsky,
{\em Mirror Symmetry for hyperk\"ahler manifolds,}
 alg-geom/9512195,  Mirror symmetry, III (Montreal, PQ, 1995), pp. 115--156,
   \textit{AMS/IP Stud. Adv. Math.}, 10, Amer. Math. Soc., Providence, RI, 1999.
   
\bibitem[Y]{_Yau:Calabi-Yau_} 
S.T. Yau, {\em On the Ricci curvature of a compact K\"ahler manifold 
and the complex Monge-Amp\`ere equation I.}  \textit{Comm. on Pure and Appl.
Math.}, \textbf{31}, pp. 339--411 (1978).


\end{thebibliography}
}
\end{document}